\documentclass[12pt]{amsart}

\usepackage{mathrsfs}
\usepackage{amsmath,amsthm,amssymb}
\font\teneufm=eufm10
\font\seveneufm=eufm7
\font\fiveeufm=eufm5
\newfam\frakturfam

\textfont\frakturfam=\teneufm
\scriptfont\frakturfam=\seveneufm
\scriptscriptfont\frakturfam=\fiveeufm

\newtheorem{lm}{Lemma}
\newtheorem{theor}{Theorem}
\newtheorem{co}{Corollary}

\def\bee{\begin{eqnarray}}
\def\bes{\begin{eqnarray*}}
\def\eee{\end{eqnarray}}
\def\ees{\end{eqnarray*}}

\def\a{\alpha}

\def\s{\sigma}
\def\t{\tau}
\def\d{\partial}
\def\Proof{{\sl Proof.}\ }

\title{Chein Automorphisms of Free Metabelian Anticommutative Algebras}

\begin{document}
\date{}
\maketitle

\begin{center}
{\bf Ruslan Nauryzbaev}\footnote{Department of Mathematics, L.N. Gumilyov Eurasian National University,
Astana, Kazakhstan,  
e-mail:{\em nauryzbaevr@gmail.com}}, 
{\bf Ivan Shestakov}\footnote{Shenzhen International Center for Mathematics, Shenzhen, China, and Instituto de Matem\'atica e Estat\'\i stica, Universidade de S\~ao Paulo, Brazil, 
e-mail:{\em shestak@ime.usp.br}}, 
and
{\bf Ualbai Umirbaev}\footnote{Department of Mathematics,
 Wayne State University,
Detroit, USA, 
and Institute of Mathematics and Mathematical Modeling, Almaty, Kazakhstan,
e-mail: {\em umirbaev@wayne.edu}}

\bigskip

\dedicatory{Dedicated to the 75th birthday of Professor V. Drensky}
\end{center}

\begin{abstract} We describe all automorphisms of a free metabelian anticommutative algebra of rank $n\geq 3$ over a field $K$ that move only one variable while fixing the others. Such automorphisms are called Chein automorphisms in the cases of free metabelian groups and free metabelian Lie algebras. 
We show that all automorphisms of a free metabelian anticommutative algebra of rank $n=2$ are linear, and that the simplest non elementary Chein automorphism of degree $3$ is absolutely  wild for all $n\geq 3$. 
\end{abstract}



\tableofcontents

\section{Introduction}

\hspace*{\parindent}

 Let ${\mathfrak M}$ be an arbitrary variety of algebras over a field $K$ and let 
 $A=K_{\mathfrak M}\langle x_1,x_2,\ldots,x_n\rangle$ be the free algebra of ${\mathfrak M}$ with a free set of generators $X=\{x_1,x_2,\ldots,x_n\}$.  As usual we denote by $n$-tuple $\phi=(f_1,\ldots,f_n)$ the endomorphism $\phi$ of $A$ such that $\phi(x_i)=f_i,\, 1\leq i\leq n$. In particular, $\mathrm{Id}=(x_1,\ldots,x_n)$ is the identity automorphism. Let $\mathrm{Aut}(A)$ be the group of all automorphisms of $A$. An automorphism  
\bee\label{f1}
(x_1,\ldots,x_{i-1}, \a x_i+f, x_{i+1},\ldots,x_n),
\eee
where $0\neq\a\in K,\ f\in K_{\mathfrak M}\langle X\setminus \{x_i\}\rangle$, 
 is called {\em elementary}.    The subgroup $\mathrm{TAut}(A)$ of $\mathrm{Aut}(A)$ generated by all 
  elementary automorphisms is called the {\em tame automorphism group}, 
 and the elements of this subgroup are called {\em tame automorphisms} 
 of $A$. Nontame automorphisms of $A$ are called {\em wild}.

An automorphism $\phi$ is called {\em approximately tame} \cite{SU25}  if it can be approximated by a sequence of 
tame automorphisms $\{\psi_k\}_{k\geq 0}$ with respect to the {\em formal power series topology}  (see exact definitions in Section 7).  An automorphism $\phi$ is called {\em absolutely wild} if it is not approximately tame.

In 1968 O. Chein \cite{Chein} introduced a specific type of automorphisms of free metabelian groups of finite rank $n\geq 3$, and proved that they are not tame when $n=3$. Nevertheless they are tame when $n\geq 4$ since all automorphisms of free metabelian groups of rank $n\geq 4$ are tame by the well known Bachmuth-Mochizuki-Roman'kov Theorem \cite{BM85,Romankov85}.  

Analogues of these automorphisms for free metabelian algebras were studied by V. Roman'kov and his coauthors  \cite{KR,Romankov08,Romankov20} as well as by R. Nauryzbaev \cite{Nauryzbaev09,Nauryzbaev11}. Let $M$ be a free metabelian Lie algebra in the variables $x_1,\ldots,x_n$.  
As usual, the group of all linear automorphisms of $M$ will be identified with the group of all invertible matrices $GL_n(K)$. In 2015, Nauryzbaev \cite{Nauryzbaev15} proved that $\mathrm{TAut}(M)$ is generated by all linear automorphisms and the quadratic automorphism 
\bee\label{f2}
(x_1+[x_2,x_3],x_2,\ldots,x_n)
\eee
if either $n\geq 4$ and $K$ is a field of characteristic $\neq 3$, or $n\geq 5$ and $K$ is a field of arbitrary characteristic.

A peculiar feature of metabelian Lie algebras is that they admit automorphisms of the form (\ref{f1}) that are not necessarily elementary; that is, $f$ may also involve $x_i$.  
We call them {\em Chein automorphisms} as in the case of groups \cite{GGR,Romankov95}. Notice that, by our definition, all elementary automorphisms form a subset of Chein automorphisms. These automorphisms move only one variable and fix the others. The subgroup $\mathrm{ATAut}(M)$ of $\mathrm{Aut}(M)$ generated by all 
  Chein automorphisms will be called the {\em almost tame automorphism group}, 
 and the elements of this subgroup will be called {\em almost tame automorphisms} 
 of $M$ \cite{Nauryzbaev09}. Thus 
\bes
\mathrm{TAut}(M)\subseteq \mathrm{ATAut}(M)\subseteq \mathrm{Aut}(M). 
\ees

The automorphisms of the algebra $M$ of rank $n=3$ have been extensively investigated by numerous authors. In 1992, V. Drensky \cite{Drensky92} (see also \cite{BD93-1}) proved that the exponential automorphism 
\bes
\mathrm{exp}(\mathrm{ad}[x_1,x_2])=(x_1+[[x_1,x_2],x_1], x_2+[[x_1,x_2],x_2], x_3+[[x_1,x_2],x_3])
\ees
 of $M$ is wild. Moreover, this automorphism is not almost tame \cite{Nauryzbaev09,Romankov08}.
Every non-elementary Chein automorphism is also shown to be wild \cite{Nauryzbaev09,Romankov08}. The defining relations for elementary automorphisms \cite{Nauryzbaev10} and for Chein automorphisms \cite{Nauryzbaev11} have also been determined. 
It was recently observed \cite{SU25} that the examples of wild automorphisms mentioned here are actually absolutely wild.

 The first example of a Chein automorphism is 
\bee\label{f3}
(x_1+[[x_2,x_3],x_1],x_2,\ldots,x_n), 
\eee
 which is absolutely wild when $n=3$ and plays an important role in the study of $\mathrm{Aut}(M)$ for $n\geq 4$.  In 1993 Bryant and Drensky \cite{BD93-2} (see also \cite{KP16}) proved that every automorphism of $M$ of rank $n\geq 4$ is approximately tame in characteristic zero. Thus (\ref{f3}) is approximately tame in characteristic zero. Recently it was shown in \cite{U24} that $\mathrm{Aut}(M_n)$ is generated modulo $GL_n(K)$ by (\ref{f2}) and  (\ref{f3}) when either $n=4$ and the characteristic of $K$ is not $3$, or $n\geq 5$ and $K$ is an arbitrary field. It remains an open question whether (\ref{f3}) is tame or wild for $n\geq 4$. 

This paper was motivated by the behavior of the Chein automorphism (\ref{f3}) in free metabelian Lie algebras. Here we study automorphisms of free metabelian anticommutative algebras. These algebras also admit Chein automorphisms. First we describe all Chein automorphisms  
of free metabelian anticommutative algebras. 
We show that every automorphism of a free metabelian anticommutative algebra of rank $2$ is tame; that is, all such automorphisms are linear. We also show that the analogue of the Chein automorphism (\ref{f3}) for free metabelian anticommutative algebras of rank $n\geq 3$ is absolutely wild.

\section{Free metabelian anticommutative algebras}

\hspace*{\parindent}

Let $K$ be an arbitrary field. An algebra $A$ over $K$ is called {\em anti-commutative} if the identity 
\bee\label{a1}
x^2=0
\eee
is satisfied. This identity implies 
\bee\label{a2}
xy=-yx. 
\eee

An algebra $A$ is called {\em metabelian} if it satisfies the identity 
\bee\label{a3}
(xy)(zt)=0. 
\eee

Let $A$ be a free metabelian anticommutative algebra in the variables $x_1,\ldots,x_n$. The identities (\ref{a1})--(\ref{a3}) imply that every element of $A^2=AA$ is a linear combination of the elements of the form 
\bee\label{a4}
(((x_i x_j) x_{k_1})\cdots )x_{k_s},
\eee
where $i>j$, $1\leq k_r\leq n$ for all $1\leq r\leq s$, and $s\geq 0$.  

Let $Y=A/A^2$ and denote by $y_i$ the images of $x_i$ in $Y$. Then $Y$ is an algebra with trivial multiplication and with a linear basis $y_1,\ldots,y_n$. Denote by $U$ the free associative algebra with identity and with free generators $z_i=R_{y_i}$, where $1\leq i\leq n$. In fact, $U$ is the universal (multiplicative) enveloping algebra of $Y$ in the variety of all anticommutative algebras \cite{U94,U96} and $z_i=R_{y_i}$ is the universal operator of right multiplication by $y_i$.  Notice that $L_{y_i}=-R_{y_i}$ by (\ref{a2}).

Let $T$ be a free right $U$-module with free generators $t_1,\dots,t_n$. Turn the space $Y \oplus T$ into an algebra by 
\bes
 (a_1+m_1)(a_2+m_2) = m_1 R_{a_2} - m_2 R_{a_1}, \quad a_1,a_2\in Y,\; m_1,m_2\in T.
\ees
Obviously, $Y \oplus T$ is a metabelian anticommutative algebra. 
\begin{lm}\label{l1} $(1)$ The homomorphism $\rho : A\to Y \oplus T$ defined by $\rho(x_i)=y_i+t_i$ for all $i$ is an embedding; 

$(2)$ The elements of the form (\ref{a4}) is a basis of $A^2$. 
\end{lm}
\Proof It is sufficient to check that the images of elements of the form (\ref{a4}) in $Y \oplus T$ under $\rho$ are linearly independent. In fact, 
\bes
 \rho(x_i x_j) = t_i z_j - t_j z_i
\ees
and 
\bes
  \rho((((x_i x_j) x_{k_1})\ldots )x_{k_s}) = (t_i z_j - t_j z_i) z_{k_1}\ldots z_{k_s}. 
\ees

Let us introduce a $\mathrm{deglex}$ order on the set of all basis elements 
\bes
t_iz_{i_1}\ldots z_{i_k}
\ees
of $T$ by setting $t_i<t_j$ and $z_i<z_j$ when $i<j$. Then the leading monomial of $\rho((((x_i x_j) x_{k_1})\ldots )x_{k_s})$ is  $t_i z_j z_{k_1}\ldots z_{k_s}$. Thus, the leading monomials of the images of all distinct elements of the form (\ref{a4}) are different and therefore linearly independent. $\Box$

By this lemma, we will identify $A$ with the subalgebra of $Y \oplus T$ generated by $x_i=y_i+t_i$, where $1\leq i\leq n$. Then every element $a\in A$ can be uniquely represented in the form 
\bes
a=y+t_1u_1+\ldots+t_nu_n, 
\ees
where $y\in Y$ and $u_1,\ldots,u_n\in U$. 

Set 
\bes
\frac{\partial a}{\partial x_i}=u_i
\ees
for all $1 \leq i \leq n$. Then $\frac{\partial a}{\partial x_i}$ are called the Fox derivatives of $a$  \cite{Shpilrain, U93,U95}. Set also 
\bes
\partial(a) = (u_1, \dots, u_n)^t = \left( \frac{\partial a}{\partial x_1}, \dots, \frac{\partial a}{\partial x_n} \right)^t, 
\ees
where $t$ denotes the transposition. We have $\partial(x_i) = e_i^t$, where  $e_i$ is the row vector with a $1$ in the $i$th place and zeros elsewhere.  

\begin{lm}\label{l2} 
$(1)$  If $g \in A^2$ then the column $\partial g = (r_1, r_2, \ldots, r_n)^t$ satisfies the conditions  
\bee\label{f8}
r_i = \sum_{j=1 }^{n} z_j w_{j}^{(i)}, \ w_{j}^{(i)} = -w_{i}^{(j)}, \ w_{i}^{(i)} = 0, 
\eee
for all $i,j$.  

$(2)$ If the column $(r_1, r_2, \ldots, r_n)^t \in U^n$ satisfies the conditions (\ref{f8}) then there exists $g \in A^2$ such that $\partial g = (r_1, r_2, \ldots, r_n)^t$.
\end{lm}
\Proof $(1)$ Since $\d$ is linear, it is sufficient to check the conditions (\ref{f8}) for an arbitrary element of the form $g = (x_s x_k) u$, where $s>k$ and $u\in U$. We have  
\bes
(x_s x_k) u = (t_s z_k-t_k z_s) u = t_s z_k u - t_k z_s u= t_s r_s + t_k r_k, 
\ees
where $r_s = z_k u = z_k w_k^{(s)}, r_k = z_s (- u) = z_s w_s^{(k)}$. Consequently, $w_k^{(s)} = - w_s^{(k)}$ and $w_{i}^{(j)}=0$ when $\{i,j\}\neq \{s,k\}$.  This proves $(1)$. 

$(2)$ Suppose that $r_1, r_2, \ldots, r_n \in U$ satisfy the conditions (\ref{f8}). Set $g= \sum_{1 \leq j < i \leq n}{(x_i x_j) w_{j}^{(i)}}$. Then 
\[
\partial (g) = \partial\left(\sum_{1 \leq j < i \leq n}{(x_i x_j) w_{j}^{(i)}} \right) = \partial\left(\sum_{1 \leq j < i \leq n}{(t_i z_j - t_j z_i) w_{j}^{(i)}} \right) =
\] 
\[
= \partial\left(\sum_{1 \leq j < i \leq n} t_i z_j w_{j}^{(i)} - \sum_{1 \leq j < i \leq n} t_j z_i w_{j}^{(i)} \right) = \partial\left(\sum_{1 \leq j < i \leq n} t_i z_j w_{j}^{(i)} - \sum_{1 \leq i < j \leq n} t_i z_j w_{i}^{(j)} \right) = 
\]
\[
= \partial\left(\sum_{i=1}^n {\sum_{j=1}^{i-1} t_i z_j w_{j}^{(i)}} - \sum_{i=1}^n {\sum_{j=i+1}^{n} t_i z_j w_{i}^{(j)}} \right) = \partial\left(\sum_{i=1}^n t_i { \left[\sum_{j=1}^{i-1} z_j w_{j}^{(i)} - \sum_{j=i+1}^{n} z_j w_{i}^{(j)}\right]} \right) =
\]
\[
= \partial\left(\sum_{i=1}^n t_i { \left[\sum_{j=1}^{i-1} z_j w_{j}^{(i)} + \sum_{j=i+1}^{n} z_j w_{j}^{(i)}\right]} \right) = \partial\left(\sum_{i=1}^n t_i \sum_{j=1}^{n} z_j w_{j}^{(i)} \right) = \partial\left(\sum_{i=1}^n t_i r_i \right) =
\]
\[
(r_1, r_2, \ldots, r_n)^t. 
\] \hfill $\Box$

Notice that if $\partial g = (r_1, r_2, \ldots, r_n)^t$ for some $g \in A^2$, then $r_n$ is uniquely determined by $r_1, r_2, \ldots r_{n-1}$. In fact, if $r_i = \sum_{j=1}^n z_j w_{j}^{(i)}$, then, by Lemma \ref{l2}, we obtain   
\bee\label{f9}
r_n = \sum_{j=1}^{n-1} (- z_j w_{n}^{(j)}).
\eee

\section{Jacobian matrices}

\hspace*{\parindent}

The Jacobian matrix of an endomorphism $\phi = (f_1, \dots, f_n)$ is defined as 
\[ J(\phi) = \left[ \frac{\partial f_j}{\partial x_i} \right]_{1 \leq i,j \leq n} = [\partial(f_1) \dots \partial(f_n)]. \]

For any $f\in A$, let $\tilde f$ denote its image in $Y=A/A^2$. Every endomorphism $\phi=(f_1,\ldots,f_n)$ of $A$ induces an  
endomorphism of $Y=A/A^2$ that sends $y_i$ to $\tilde f_i$ for all $i$. This endomorphism can be uniquely extended to an endomorphism of $U$ that sends $z_i=R_{y_i}$ to $R_{\tilde f_i}$ for all $i$. Denote this endomorphism of $U$ by $\tilde \phi$. For  $f\in U$, the image  of $f$ under $\tilde \phi$ will be written as $f^{\tilde \phi}$.

If  $g$ is an arbitrary element of $A$, then  it is easy to check that the following form of the "chain rule" holds: 
\[
\partial(\phi(g))=\partial(g(f_1,\ldots, f_n)) = \sum_{1 \leq i \leq n} \partial(f_i)  (\frac{\partial g}{\partial x_i})^{\tilde \phi}. 
\]
Recall that the composition of endomorphisms $\phi,\psi$ of $A$ is defined as  
\[
\phi \circ \psi (a) = \phi(\psi(a)), \ a\in A. 
\]

The chain rule immediately implies that 
\bee\label{f10}
J(\phi \circ \psi) = J(\phi) (J(\psi))^{\tilde \phi}. 
\eee

Consequently, if $\phi$ is an automorphism of $A$, then $J(\phi)$ is invertible over $U$. The converse of this statement is an analogue of the famous Jacobian Conjecture \cite{Keller} for free metabelian anticommutative algebras \cite{U94,U96}. 
Following the arguments in \cite{U93,U95},  one can show that the following lemma holds. 
\begin{lm}\label{l3} 
An endomorphism $\phi$ of $A$ is an automorphism if and only if the Jacobian matrix $J(\phi)$ is invertible.  
\end{lm}
We omit the proof of this lemma, as it is not needed in the present paper.

Denote by  $\mbox{IAut}(A)$ the group of all automorphisms of $A$ that induces the identity automorpjism on $Y=A/A^2$. 
The elements of $\mbox{IAut}(A)$ are called {\em $\mbox{IA}$-automorphisms}. If $\phi\in \mbox{IAut}(A)$ then the Jacobian matrix of $\phi$ has the form 
\[
J(\phi) = I_n+R,
\]
where $I_n$ is the identity matrix of order $n$ and every column of $R$ satisfies the conditions (\ref{f8}). If $\phi\in\mbox{IAut}(A)$ then $\tilde \phi=\mbox{Id}$. Consequently, if $\phi,\psi\in\mbox{IAut}(A)$ then (\ref{f10}) implies that 
\bes
J(\phi \circ \psi) = J(\phi) J(\psi). 
\ees
\begin{co}\label{c1} The mapping 
\bes
J: \mbox{IAut}(A)\to SL_n(U)
\ees
is an embedding of groups. 
\end{co}
It is obvious that $J$ is one-to-one. Recall that $SL_n(U)$ denotes the set of all invertible matrices over $U$ whose constant part has determinant one.

\section{The structure of Chein automorphisms}

\hspace*{\parindent}

\begin{lm}\label{l4} 
Let $f\in A^2$. Then $\frac{\d f}{\d x_1}=0$ if and only if $f$ belongs to the ideal of $A$ generated by all $x_ix_j$, where $i,j>1$. 
\end{lm}
Let $f\in A^2$. Then 
\bes
f= \sum_{1\leq j<i \leq n} (x_i x_j) u_{ij}, \ u_{ij}\in U. 
\ees
Consequently, 
\bes
\frac{\d f}{\d x_1}= -\sum_{i=2}^nz_iu_{i1}. 
\ees
If $\frac{\d f}{\d x_1}=0$ in the free associative algebra, we immediately get $u_{21}=\ldots =u_{n1}=0$. Therefore 
\bee\label{f11}
f= \sum_{2\leq j<i \leq n} (x_i x_j) u_{ij}, 
\eee
satisfies the statement of the lemma. $\Box$

Consider the endomorphisms of the form  
\bee\label{f12}
\delta=(x_1+f,x_2,\ldots,x_n), \ f\in A^2. 
\eee
\begin{lm}\label{l5} 
An endomorphism $\delta$ of $A$ of the form (\ref{f12}) is an automorphism if and only if $\frac{\d f}{\d x_1}=0$. 
\end{lm}
\Proof We have 
\bes
J(\delta)=\begin{bmatrix}
 1 + \frac{\partial f}{\partial x_1} & 0 & \ldots & 0 \\
\frac{\partial f}{\partial x_2} & 1 & \ldots & 0 \\
\cdot & \cdot & \ldots & \cdot \\
 \frac{\partial f}{\partial x_n} & 0 & \ldots & 1 
\end{bmatrix}. 
\ees
If $J(\delta)$ is invertible then $1 + \frac{\partial f}{\partial x_1}$ is right invertible in the free associative algebra $U$. Since $\frac{\partial f}{\partial x_1}$ belongs to the augmented ideal of $U$ it follows that $\frac{\partial f}{\partial x_1}=0$. By Lemma \ref{l4}, $f$ has the form  (\ref{f11}). Then $\delta$ is a product of endomorphisms 
\bes
\delta_{ij}=(x_1+(x_i x_j) u_{ij},x_2,\ldots,x_n),  \ 2\leq j<i \leq n. 
\ees
Since $A$ is metabelian it follows that every $\delta_{ij}$ is an automorphism. $\Box$

\section{Automorphisms of algebras of rank $2$}

\hspace*{\parindent}

If $n=1$ then the group of automorphisms  $\mbox{Aut}(A)$ is isomorphic to $GL_1(K)\simeq K^*=K\setminus \{0\}$. If $n=2$ then all elementary automorphisms are linear since $A$ is anticommutative. Consequently, $\mbox{TAut}(A)\simeq GL_2(K)$. We show that $A$ has no other automorphisms.  
\begin{lm}\label{l6} 
If $n=2$ then $\mbox{Aut}(A)\simeq GL_2(K)$. 
\end{lm}
\Proof It is sufficient to prove that $\mbox{IAut}(A)=\{\mbox{Id}\}$. Let $\phi=(x_1+f_1,x_2+f_2) \in \mbox{IAut} (A)$, where $f_1,f_2\in A^2$. Then   
\[
J(\varphi) = I_2 +
\begin{bmatrix}
\d(f_1) & \d(f_2) 
\end{bmatrix} = I_2 +
\begin{bmatrix}
a & b \\
c & d
\end{bmatrix}.
\]
By Lemma \ref{l2}, we have 
\[
a = z_2 h, \ c = -z_1 h, \ h\in U.  
\]
The first column of $J(\varphi)$ is $[1+ z_2 h, -z_1 h]^t$. Since $J(\varphi)$ is invertible it follows that $[1+ z_2 h, -z_1 h]^t$ is left unimodular, that is, there exist $u_1,u_2\in U$ such that 
\bee\label{f13}
u_1(1+ z_2 h)+u_2(-z_1 h)=1. 
\eee
For any $u\in U$ denote by $\overline{u}$ its highest homogeneous part of $u$ in $U$ with respect to the degree function in the variables $z_1,z_2$. If $h\neq 0$, then $\overline{1+ z_2 h}=z_2\overline{h}$ and $\overline{z_1 h}=z_1\overline{h}$. Then (\ref{f13}) implies that 
\bes
\overline{u_1}z_2\overline{h}-\overline{u_2}z_1\overline{h}=0, 
\ees
that is, $z_2\overline{h}$ and  $z_1\overline{h}$ are left dependent. It is well known that in free associative algebras this implies \cite{Cohn06} that one of $z_2\overline{h}$ or $z_1\overline{h}$ belongs to the left ideal generated by the other. They must be linearly dependent since they have the same degrees. But they are linearly independent if $h\neq 0$. 

Consequently, $h=0$ and $a=c=0$. Similarly, $b=d=0$. Then $f_1=f_2=0$ and $\phi=\mbox{Id}$. 
$\Box$

\section{Divergence in free anticommutative algebras}

\hspace*{\parindent}

To prove the main result of the paper, we need the notion of the divergence of a derivation introduced in \cite{SU25}.

Let $B$ be a free anticommutative algebra over $K$ in the variables $\Xi=\{\xi_1,\ldots,\xi_n\}$. Let $\Xi^*$ be the set of all (nonempty) nonassociative monomials in the alphabet $\Xi$. Recall that every $u\in \Xi^*$ of degree $\geq 2$ can be uniquely written as $u=u_1u_2$ for some $u_1,u_2\in \Xi^*$ \cite{KBKA}. Set $\xi_i<\xi_j$ if $i<j$. For any $u,v\in \Xi^*$ set $u<v$ if $\deg u<\deg v$. If $\deg u=\deg v\geq 2$ and $u=u_1u_2$ and $v=v_1v_2$, then set $u<v$ if either $u_1<v_1$ or $u_1=v_1$ and $u_2<v_2$. 

A nonassociative monomial $u\in \Xi^*$ is called {\em regular} if it does not have submonomials of the form $vw$, where $v,w\in  \Xi^*$ and $v\leq w$. The set of all regular words is a linear basis of $B$ \cite{Shirshov54}

The universal (multiplicative) enveloping algebra $U(B)=U_{\mathfrak{M}}(B)$ is the free associative algebra with identity generated by all universal operators of right multiplication $R_u$, where $u$ runs over the set of all regular monomials in the alphabet $\Xi$ \cite{U94}.

Let 
\bes
\Omega_{B}=d\xi_1U(B)\oplus \ldots \oplus d\xi_nU(B)
\ees
be a free right $U(B)$-module with free generators $d\xi_1,\ldots,d\xi_n$. 
 The linear map 
\bes
D : B \to \Omega_{B}
\ees
defined by $D(\xi_i)=d\xi_i$ for all $i$ and such that 
\bes
D(ab)=aD(b)+D(a)b= -D(b)R_a+ D(a)R_b
\ees
for all $a,b\in B$, is called the {\em universal derivation} of $B$ and $\Omega_{B}$ is called the {\em universal differential module} of $B$ \cite{U94,U96}. 

For any $a\in B$ there exist unique elements $u_1,u_2,\ldots,u_n\in U(B)$ such that
\bes
D(a)=d\xi_1u_1+\ldots+d\xi_nu_n.
\ees
The elements $u_i=\frac{\partial a}{\partial \xi_i}$ are called the {\em Fox derivatives} of $a\in B$ \cite{U94,U96}. Set also 
\bes
\d(a)=(\frac{\partial a}{\partial \xi_1},\ldots,\frac{\partial a}{\partial \xi_n})^t. 
\ees

Consider the natural grading 
\bes
 B= B_1\oplus \ldots \oplus B_k\oplus\ldots 
\ees
of $B$ with respect to the standard function $\deg$, that is, $B_k$ is the linear span of (non-associative) monomials of degree $k$. 

Let $L=\mathrm{Der} (B)$ be the Lie algebra of all derivations of $B$. For any $n$-tuple $F=(f_1,\ldots,f_n)\in B^n$ there exists a unique derivation $D$ of $B$ such that $D(\xi_i)=f_i$ for all $1\leq i\leq n$. Denote this derivation by 
\bee\label{f14}
D=f_1\partial_1+\ldots+f_n\partial_n=D_F. 
\eee
For any derivation of this form define its Jacobian matrix by  
\bes
J(D)=[\frac{\partial f_j}{\partial \xi_i}]_{1\leq i,j\leq n}= [\partial(f_1) \dots \partial(f_n)].
\ees

We say that $D$ from (\ref{f14}) is homogeneous of degree $i$ if $f_1,\ldots,f_n\in B_{i+1}$. Let $L_i$ be the space of all homogeneous derivations of degree $i$. Then 
\bes
L=L_0\oplus L_1\oplus\ldots\oplus L_k\oplus \ldots, \ \ \ [L_i,L_j]\subseteq L_{i+j}, 
\ees
is a grading of $L$. Notice that $L_0$ is isomorphic to the matrix algebra $\mathfrak{gl}(n)$ with the matrix units $e_{ij}=\xi_i\partial_j$, where $1\leq i,j\leq n$. 

The {\em divergence} $\mathrm{div}(D)$ of $D$ is defined as the image of the trace of the Jacobian matrix $J(D)$ in the quotient space $U(B)/[U(B),U(B)]$ \cite{SU25}. This agrees with the definition of the divergence given in \cite[Definition 1]{SU25}, since the Jacobson radical $R=R(U)$ of the free associative algebra $U$ is zero. It is shown in \cite{SU25} that 
the space 
\bes
S=S_0\oplus S_1\oplus\ldots\oplus S_k\oplus \ldots,
\ees
of all derivations of $B$ with zero divergence is a graded subalgebra of $L$ and called the {\em special Lie algebra of derivations} of $B$.

 Let $\mathrm{IE}(k)=\mathrm{IE}(B,k)$ be the 
 monoid of all endomorphisms of $B$ that induces the identity automorphism on the factor-algebra $B/B^{k+1}$. 
Then we get the descending series of monoids 
\bes
\mathrm{IE}(1)\supseteq\mathrm{IE}(2)\supseteq\ldots\supseteq \mathrm{IE}(k)\supseteq . 
\ees
If $\phi\in \mathrm{IE}(i)\setminus\mathrm{IE}(i+1)$ for some $i\geq 1$, then  
\bes
\phi=(\xi_1+f_1+F_1,\ldots,\xi_n+f_n+F_n), 
\ees
where $f_1,\ldots,f_n$ are homogeneous elements of degree $i+1$, with at least one of them nonzero,  and $F_j\in B^{i+2}$ for all $1\leq j\leq n$. Set 
\bes
T(\phi)=f_1\d_1+\ldots+f_n\d_n. 
\ees
The derivation $T(\phi)$ is referred to as the tangent derivation of the endomorphism $\phi$ (with respect to the power series topology) \cite{SU25}.

Consider the endomorphism 
\bee\label{f15}
\s=(\xi_1+(\xi_3\xi_2)\xi_1,\xi_2,\ldots,\xi_n)
\eee
\begin{lm}\label{l7} 
$\mathrm{div}(T(\s))\neq 0$. 
\end{lm}
\Proof We have $T(\s)=((\xi_3\xi_2)\xi_1)\d_1$. Then the first entry on the main diagonal of $J(T(\s))$ is 
\bes
\frac{\partial (\xi_3\xi_2)\xi_1)}{\partial \xi_1}, 
\ees
and all the other diagonal entries are zero. Consequently, the trace of the matrix $J(T(\s))$ is 
\bes
\frac{\partial (\xi_3\xi_2)\xi_1}{\partial \xi_1}=-R_{\xi_3\xi_2}\notin [U(B),U(B)]. 
\ees
Therefore $\mathrm{div}(T(\s))\neq 0$. $\Box$

\section{An absolutely wild automorphism}

\hspace*{\parindent}

In this section we show that the Chein automorphism 
\bee\label{f16}
\t=(x_1+(x_3x_2)x_1,x_2,\ldots,x_n)
\eee
of the free metabelian anticommutative algebra $A$ is absolutely wild for all $n\geq 3$. 

First, we recall the precise definitions of approximately tame and absolutely wild automorphisms from \cite{SU25}.

 Let $\mathrm{IA}(k)=\mathrm{IA}(A,k)$ be the 
 set of all automorphisms of $A$ that induces the identity automorphism on the factor-algebra $A/A^{k+1}$. 
Consider the descending central series 
\bes
\mathrm{IAut}(A)=\mathrm{IA}(1)\supseteq\mathrm{IA}(2)\supseteq\ldots\supseteq \mathrm{IA}(k)\supseteq . 
\ees
 An automorphism $\phi\in \mathrm{Aut}(A)$ is called {\em approximately tame} if there exists a sequence of 
tame automorphisms $\{\psi_k\}_{k\geq 0}$ such that $\phi\psi_k^{-1}\in \mathrm{IA}(k)$. The topology defined by this definition on $\mathrm{Aut}(A)$ is called the {\em formal power series topology}  \cite{Anick}. An automorphism $\phi$ will be called {\em absolutely wild} if it is not approximately tame.  

In 1981 Shafarevich \cite{Shafarevich81} and in 1983 Anick \cite{Anick} independently proved that every automorphism of the polynomial algebra $K[x_1,x_2,\ldots,x_n]$ over a field $K$ of characteristic zero is approximately tame. 
In 1993 Bryant and Drensky \cite{BD93-2} (see also \cite{KP16}) proved that every automorphism of $M$ of rank $n\geq 4$ is approximately tame in characteristic zero. 

The following result is a part of Theorem 5 from \cite{SU25}. 

$(F)$ Let $\mathfrak{M}$ be  a Nielsen-Schreier variety of algebras and $B$ is a free algebra of $\mathfrak{M}$ of rank $n$. 
Let $\mathfrak{N}$ be any subvariety of $\mathfrak{M}$ and let $I\subseteq B$ be the ideal of all identities of $\mathfrak{N}$ in $\mathfrak{M}$ in $n$ variables. Suppose that an endomorphism $\epsilon\in \mathrm{IE}_i(A)\setminus\mathrm{IE}_{i+1}(A)$ for some $i\geq 1$ induces an automorphism $\phi$ of $B=A/I$ and the ideal $I$ does not contain elements of degree $\leq i+1$.
 If $\mathrm{div}(T(\epsilon))\neq 0$ then $\phi$ is absolutely wild. 

\begin{theor}\label{t1}
The automorphism $\t$ from (\ref{f16}) is an absolutely wild automorphism of the free metabelian anticommutative algebra $A$ of rank $n\geq 3$. 
\end{theor}
\Proof We apply $(F)$. Let $\mathfrak{M}$ be the variety of all anticommutative algebras. It is well known that $\mathfrak{M}$ is Nielsen-Schreier \cite{Shirshov54}. Let $B$ be the free anticommutative algebra in the variables $\Xi=\{\xi_1,\ldots,\xi_n\}$. Let $\mathfrak{N}$ be the subvariety of $\mathfrak{M}$ defined by the identity (\ref{a3}). Then $I$ is the ideal of $B$ generated by all elements of the form $(uv)(wq)$, where $u,v,w,q\in \Xi^*$. 

Notice that the endomorphism $\s$ of $B$ induces the automorphism $\t$ of $A=B/I$. We have $\s\in \mathrm{IE}_2(A)\setminus\mathrm{IE}_3(A)$. The ideal $I$ does not contain elements of degree $\leq 3$. Moreover, $\mathrm{div}(T(\s))\neq 0$ by Lemma \ref{l7}. By $(F)$, $\t$ is an absolutely wild automorphism of $A$. 

Note that, in general, all results in \cite{SU25} are stated for fields of characteristic zero. By carefully repeating the proof of Theorem 5 from \cite{SU25} for these specific varieties $\mathfrak{M}$ and $\mathfrak{N}$, we see that the assumption of characteristic zero is not required. Thus, $\t$ is an absolutely wild automorphism of $A$ of rank $n\geq 3$ over an arbitrary field.  $\Box$

\section*{Acknowledgments}

The first and the third authors are supported by grant AP23486782 from the Ministry of Science and Higher Education of the Republic of Kazakhstan.

The second author is supported by ICM, SUSTech, Shenzhen, China, and by Brazilian grants FAPESP 2024/14914-9 and CNPq 305196/2024-3.

\end{document}